\PassOptionsToPackage{unicode}{hyperref}
\PassOptionsToPackage{hyphens}{url}
\documentclass[
]{article}
\usepackage{amsmath,amssymb}
\usepackage{lmodern}
\usepackage{iftex}
\ifPDFTeX
  \usepackage[T1]{fontenc}
  \usepackage[utf8]{inputenc}
  \usepackage{textcomp} 
\else 
  \usepackage{unicode-math}
  \defaultfontfeatures{Scale=MatchLowercase}
  \defaultfontfeatures[\rmfamily]{Ligatures=TeX,Scale=1}
\fi
\IfFileExists{upquote.sty}{\usepackage{upquote}}{}
\IfFileExists{microtype.sty}{
  \usepackage[]{microtype}
  \UseMicrotypeSet[protrusion]{basicmath} 
}{}
\makeatletter
\@ifundefined{KOMAClassName}{
  \IfFileExists{parskip.sty}{%
    \usepackage{parskip}
  }{
    \setlength{\parindent}{0pt}
    \setlength{\parskip}{6pt plus 2pt minus 1pt}}
}{
  \KOMAoptions{parskip=half}}
\makeatother
\usepackage{xcolor}
\usepackage{longtable,booktabs,array}
\usepackage{multirow}
\usepackage{calc} 
\usepackage{etoolbox}
\makeatletter
\patchcmd\longtable{\par}{\if@noskipsec\mbox{}\fi\par}{}{}
\makeatother
\IfFileExists{footnotehyper.sty}{\usepackage{footnotehyper}}{\usepackage{footnote}}
\makesavenoteenv{longtable}
\usepackage{graphicx}
\makeatletter
\def\maxwidth{\ifdim\Gin@nat@width>\linewidth\linewidth\else\Gin@nat@width\fi}
\def\maxheight{\ifdim\Gin@nat@height>\textheight\textheight\else\Gin@nat@height\fi}
\makeatother
\setkeys{Gin}{width=\maxwidth,height=\maxheight,keepaspectratio}
\makeatletter
\def\fps@figure{htbp}
\makeatother
\ifLuaTeX
  \usepackage{selnolig}  
\fi
\IfFileExists{bookmark.sty}{\usepackage{bookmark}}{\usepackage{hyperref}}
\IfFileExists{xurl.sty}{\usepackage{xurl}}{} 
\hypersetup{
  hidelinks,
  pdfcreator={LaTeX via pandoc}}

\author{}
\date{}

\begin{document}

UDC 517.956.35

{\large \textbf{ON THE ABSENCE AND NON-UNIQUENESS OF CLASSICAL SOLUTIONS OF MIXED PROBLEMS FOR THE TELEGRAPH EQUATION WITH A NONLINEAR POTENTIAL}}

\textbf{V. I. Korzyuk\textsuperscript{1}, J. V. Rudzko\textsuperscript{2}}

\emph{\textsuperscript{1}Institute of Mathematics of the National Academy of Sciences of Belarus, Belarus, Minsk,} \texttt{korzyuk@bsu.by}

\emph{\textsuperscript{2}Institute of Mathematics of the National Academy of Sciences of Belarus, Belarus, Minsk,} \texttt{janycz@yahoo.com}

For the telegraph equation with a nonlinear potential given in the first
quadrant, we consider the first and the second mixed problem, for which
we study issues related to the absence and non-uniqueness of classical solutions.

\emph{\textbf{Keywords:}} semilinear wave equation, mixed problem,
classical solution, absence of solution, non-uniqueness of solution, method of characteristics, energy methods, matching conditions.

\textbf{1. Introduction}

Continuous media are described mainly by nonlinear partial differential
equations. The choice of linear or nonlinear equations for describing a
medium depends on the role played by nonlinear effects and is determined
by the specific physical situation. For example, when describing the
propagation of laser pulses, it is necessary to take into account the
dependence of the refractive index of the medium on the electromagnetic field intensity. The linearization of nonlinear equations of mathematical physics does not always lead to meaningful results. It may turn out that the linearized equations apply to the physical process in question only for some finite time. Moreover, from the viewpoint of physics, it is often ``essentially nonlinear'' solutions, qualitatively different from the solutions of linear equations, that are extremely important for nonlinear equations of mathematical physics. These can be stationary solutions of the soliton type, localized in one or several dimensions, or solutions of the wave collapse type, which describe the spontaneous concentration of energy in small regions of space. Stationary solutions of hydrodynamic equations are also essentially nonlinear~[1].

Nonlinear partial differential equations are difficult to study: almost
no general techniques exist that work for all such equations, and
usually each individual equation has to be studied as a separate
problem. A fundamental question for any partial differential equation
is the existence and uniqueness of a solution for given boundary
conditions. For nonlinear equations these questions are in general very
hard~[2].

This paper deals with the question of the absence and non-uniqueness of global and local classical solutions of the telegraph equation with a nonlinear potential.

\textbf{2. Statement of the problem}

In the domain $Q=(0,\infty )\times (0,\infty )$ of two independent variables $(t,x)\in (0,\infty )\times (0,\infty )\subset {{\mathbb{R}}^{2}}$, consider the one-dimensional nonlinear equation.
\begin{equation}
\left( \partial _{t}^{2}-{{a}^{2}}\partial _{x}^{2} \right)\!u\!\left( t,x \right)-f\!\left( t,x,u\!\left( t,x \right) \right)=F\!\left( t,x \right),\quad \left( t,x \right)\in Q,
\end{equation}
where $a\in (0,\infty )$, $F$ is a function given on the set $\overline{Q}$, and $f$ is a function given on the set $\overline{Q} \times \mathbb{R}$. Eq.~(1) is equipped with the initial condition
\begin{equation}
u(0,x)=\varphi (x),\quad {{\partial }_{t}}u(0,x)=\psi (x),\quad x\in [0,\infty ),
\end{equation}
and the boundary condition
\begin{equation}
B[u](t,0)=\mu (t),\quad t\in [0,\infty ).
\end{equation}
where $\varphi$, $\psi$, and $\mu$ are functions given on the half-line $[0,\infty)$ and $B$ is some operator (it can have various forms, but in the present paper we assume that $B = I$ or $B = \partial_x$, where $I$ is the identity operator)

\textbf{3. Nonexistence of solutions for inhomogeneous matching
conditions}

\textbf{Assertion 1.} {\it Assume $B=I$. If the homogeneous matching conditions
\begin{gather*}
\mu (0)=\varphi (0),\quad {\mu }'(0)=\psi (0), \\
\mu''(0)=\frac{1}{2}\Big( f\left( 0,0,\varphi (0) \right)+f\left(0,0,\mu(0)  \right) \!\Big)+F(0,0)+{{a}^{2}}\varphi (0)
\end{gather*}
fail for given functions $f$, $\varphi$, $\psi$, $\mu$, and $F$, then the first mixed problem (1) -- (3) does not have a classical solution defined on $\overline{Q}$. }

The \textbf{proof} follows from Theorem 1 from~[3, 4].

\textbf{Assertion 2.} {\it Assume $B=\partial_x$. If the homogeneous matching conditions
$$
\mu(0)=\varphi'(0),\quad\mu'(0)=\psi'(0),
$$
fail for given functions $f$, $\varphi$, $\psi$, $\mu$, and $F$, then, for any
smoothness of these functions the second mixed problem (1) -- (3) does not have a classical solution defined on $\overline{Q}$. }

The \textbf{proof} can be carried out by the method of characteristics
by analogy with that of the preceding assertion.

\textbf{4. Nonexistence of solutions for negative energy}

In this section, we impose the following restrictions on the nonlinearity, the right side of the equation, the initial data, and the boundary data of the problems.

\textbf{Condition 1.} The functions $F$ and $\mu$ are identically equal
to zero, the function $f$ has the form $f(t,x,z)=-g(z)$, where $g(0)=0$, and the smoothness conditions $\varphi \in C_{c}^{2}([0,\infty ))$, $\psi \in C_{c}^{2}([0,\infty ))$, $g\in {{C}^{1}}(\mathbb{R})$ are satisfied.

Under condition 1, we introduce the notation
$$
G(z)=\int\limits_{0}^{z}{g(\xi )d\xi },\quad z\in \mathbb{R},
$$
and define the \emph{energy} of a solution $u$ of the problem (1) -- (3)
\begin{multline*}
E:\left[ 0,\infty  \right)\ni t\mapsto E\left( t \right)= \\
=\int\limits_{0}^{\infty }{\left( \frac{1}{2}\left( {{\left( \frac{\partial u}{\partial t}\left( t,x \right) \right)}^{2}}+{{a}^{2}}{{\left( \frac{\partial u}{\partial x}\left( t,x \right) \right)}^{2}} \right)+G\left( u\left( t,x \right) \right) \right)dx}\in \mathbb{R}.
\end{multline*}

\textbf{Assertion 3.} {\it Assume that Condition 1 holds, for some constant the inequality
$$
z g(z) \leqslant \lambda G(z),\quad z\in \mathbb{R}.
$$
Suppose also that the energy is negative
$$
E\left( 0 \right)=\int\limits_{0}^{\infty }{\left( \frac{1}{2}\left( {{\left( {\varphi }'\left( x \right) \right)}^{2}}+{{a}^{2}}{{\left( \psi \left( x \right) \right)}^{2}} \right)+G\left( \varphi \left( x \right) \right) \right)dx}<0.
$$
Then the problem (1) -- (3) does not have a classical solution defined on $\overline{Q}$. }

The \textbf{proof} of the assertion can be carried out by the scheme set
forth in~[5, Ch.~12].

\textbf{5. Nonuniqueness of solutions for nonlinearity satisfying the
H\"{o}lder condition}

In this section, we consider the second mixed problem (1) -- (3) in the following case
\begin{equation}
\begin{gathered}
f(t,x,z) := z^{\alpha}, \quad 0< \alpha < 1, \\
F(t, x) = 0,\quad \varphi = \psi = \mu \equiv 0, \quad B = \partial_x.
\end{gathered}
\end{equation}

It is easy to see that the mixed problem (1) -- (4) has the trivial
solution $u\equiv 0$. To find non-trivial solutions of the problem (1) -- (4) consider the ansatz
\begin{equation}
u(t,x)=u(t)=\beta {{t}^{\gamma }},\quad (t,x)\in \overline{Q}.
\end{equation}
where $\beta$ and $\gamma$ are some real numbers. Substituting ansatz~(5) into Eq.~(1), we obtain the relation
$$
\beta (\gamma -1)\gamma {{t}^{\gamma -2}}={{\beta }^{\alpha }}{{t}^{\alpha \gamma }},
$$
which leads to the system of equations
$$\gamma -2=\gamma \alpha ,\quad \beta (\gamma -1)\gamma ={{\beta }^{\alpha }},$$
which has the solution
\begin{equation}
\beta ={{2}^{\tfrac{1}{\alpha -1}}}{{\left( \alpha -3+\frac{4}{\alpha +1} \right)}^{\tfrac{1}{1-\alpha }}},\quad \gamma =\frac{2}{1-\alpha }.
\end{equation}

Substituting (6) into (5), we get the function
\begin{equation}
{{u}_{p}}(t,x)={{2}^{\tfrac{1}{\alpha -1}}}{{\left( \frac{\alpha +1}{{{\alpha }^{2}}-2\alpha +1} \right)}^{\tfrac{1}{\alpha -1}}}{{t}^{\tfrac{2}{1-\alpha }}}.\end{equation}

It is easy to see that the function (7) satisfies the initial (2) and
boundary conditions (3) of the problem (1)~--~(4). Thus, we have constructed one nontrivial solution of the problem (1)~--~(4), which is determined by the formula
(7). Moreover, it can be easily shown that the `glued'~[6] solution
$${{u}_{p;s}}(t,x)=\left\{ \begin{array}{*{35}{l}}
   0, & t\in [0,s),  \\
   {{u}_{p}}(t-s,x), & t\in [s,+\infty ),  \\
\end{array} \right.$$
with parameter $s>0$ also satisfies the problem (1)~--~(4). Thus, we have constructed the infinite set of nontrivial classical solutions of the problem (1)~--~(4).

We note that in the problem (1)~--~(4) the nonlinearity $u\mapsto -{{u}^{\alpha }}$ is not differentiable function on the set $\mathbb{R}$. It is the fact that makes the construction of a unique local classical solution of the problem (1)~--~(4) impossible because, in the case of continuously differentiable nonlinearity, we can build a local classical solution (but the matching conditions have to be satisfied). We can do this using the methods proposed in the works~[3--5, 7, 8].

We can use this approach to prove the non-uniqueness of the classical solution of other problems~[9, 10].

\textbf{6. Conclusions}

In this paper, we have shown that the fulfillment of the smoothness
conditions and the matching conditions is not enough for the existence
of a global classical solution of boundary value problems for the
telegraph equation with a nonlinear potential, unlike the linear
telegraph equation~[11--17]. Also, these conditions do not lead to the uniqueness of solutions. For the existence and uniqueness of a global classical solution, the nonlinearity of the equation has to satisfy some additional conditions, e.g., the Lipschitz condition~[3--5, 8, 18]. But on the other hand, the Lipschitz condition is not necessary for the existence of a unique global classical solution~[19].

\textbf{References}

\begin{enumerate}
\def\labelenumi{\arabic{enumi})}
\item Prokhorov A. M. [et al.], eds. \emph{Encyclopedia of Physics: in 5
vol}. Moscow, 1992, vol. 3. 642 p. (in~Russian).

\item \emph{Nonlinear Partial Differential Equation}. Wikipedia. Available at:
  \url{https://en.wikipedia.org/wiki/Nonlinear\_partial\_differential\_equation}
  (accessed 17 March 2023).

\item Korzyuk V. I., Rudzko J. V. Classical Solution of the First Mixed
Problem for the Telegraph Equation with a Nonlinear Potential.
\emph{Differential Equations}, 2022, vol. 58, no.~2, pp. 175--186.

\item Korzyuk V. I., Rudzko J. V. Classical and Mild Solution of the First
Mixed Problem for the Telegraph Equation with a Nonlinear Potential.
\emph{The Bulletin of Irkutsk State University. Series Mathematics},
2023, vol.~43, no.~1, pp. 48--63.

\item Evans L. C. \emph{Partial Differential Equations}. 2\textsuperscript{nd} ed. Providence, Am. Math. Soc., 2010.

\item Amel'kin V. V. \emph{Differential Equations}. Minsk, BSU, 2012. 288 p. (in~Russian).

\item
  Korzyuk V., Rudzko J. Local Classical Solution of the Cauchy problem
  for a Semilinear Hyperbolic Equation in the Case of Two Independent
  Variables. \emph{International Scientific Conference ``Ufa Autumn
  Mathematics School -- 2022'': Proceedings of the International
  Scientific Conference, Ufa, September 28 -- October 01, 2022}. Ufa,
  2022, vol. 2, pp.~48--50. (in~Russian).
  
\item
  Korzyuk V. I., Rudzko J. V. Classical Solution of the Initial-Value
  Problem for a One-Dimensional Quasilinear Wave Equation. \emph{Doklady
  of the National Academy of Sciences of Belarus}, 2023, vol. 67, no. 1,
  pp. 14--19.
  
\item
  Jokhadze O. M. Global Cauchy Problem for Wave Equations with a
  Nonlinear Damping Term. \emph{Differential Equations}, 2014, vol.~50,
  no. 1, pp.~57--65.
\item
  Jokhadze O. M. On Existence and Nonexistence of Global Solutions of
  Cauchy--Goursat Problem for Nonlinear Wave Equation. \emph{J. Math.
  Anal. Appl.}, 2008, vol. 340, no.~2, pp. 1033--1045.

\item Korzyuk V. I., Rudzko J. V. Method of Reflections for the
Klein--Gordon Equation. \emph{Doklady of the National Academy of
Sciences of Belarus}, 2022, vol. 66, no. 3, pp. 263--268.

\item Lomautsau F. E. The First Mixed Problem for the General Telegraph
Equation with Variable Coefficients on the Half-Line. \emph{Journal of
the Belarusian State University. Mathematics and Informatics}. 2021,
vol. 1, pp. 18--38 (in~Russian).

\item Lomovtsev F. E. The Second Mixed Problem for the General Telegraph
Equation with Variable Coefficients in the First Quarter of the Plane.
\emph{Vesnik of Yanka Kupala State University of Grodno. Series 2. Mathematics. Physics. Informatics, Computer Technology and Control}, 2022, vol.~12, no.~3, pp.~50--70 (in~Russian).

\item Lomovtsev F. E. Global Correctness Theorem to the First Mixed Problem
for the General Telegraph Equation with Variable Coefficients on a
Segment. \emph{Problems of Physics, Mathematics and Technics}, 2022,
vol.~50, no.~1, pp.~62--73.

\item Lomovtsev F. E., Spesivtseva K. A. Mixed Problem for a General 1D
Wave Equation with Characteristic Second Derivatives in a Nonstationary
Boundary Mode. \emph{Math. Notes}, 2021, vol. 110, no. 3, pp. 329--338.

\item Korzyuk V. I., Stolyarchuk I. I. Classical Solution of the First
Mixed Problem for the Klein--Gordon--Fock Equation in a Half-Strip.
\emph{Differential Equations}, 2014, vol. 50, no. 8, pp. 1098--1111.

\item Korzyuk V. I., Stolyarchuk I. I. Classical Solution of the First
Mixed Problem for Second-Order Hyperbolic Equation in Curvilinear
Half-Strip with Variable Coefficients. \emph{Differential Equations},
2017, vol. 53, no. 1, pp. 74--85.

\item Korzyuk V. I., Kovnatskaya O. A., Sevastyuk V. A. Goursat's Problem
on the Plane for a Quasilinear Hyperbolic Equation. \emph{Doklady of the
National Academy of Sciences of Belarus}, 2022, vol. 66, no. 4, pp.
391--396 (in~Russian).

\item
  Kharibegashvili S. S., Jokhadze O. M. Cauchy Problem for a Generalized Nonlinear Liouville Equation. \emph{Differential Equations}, 2011, vol.~47, no. 12, pp.~1741--1753.
  
\end{enumerate}

\end{document}